\documentclass[11pt,a4paper]{amsart}
\usepackage[T1]{fontenc}
\usepackage[utf8]{inputenc}
\usepackage[english]{babel}
\usepackage[a4paper,top=2.2cm,bottom=2.2cm,left=2.6cm,right=2.6cm]{geometry}
\usepackage{amsmath,amssymb,amsthm,mathtools}
\usepackage{mathrsfs}
\usepackage[numbers,sort&compress]{natbib}
\usepackage[colorlinks=true,linkcolor=black,citecolor=black,urlcolor=black]{hyperref}
\usepackage{enumitem}
\usepackage{booktabs}
\usepackage{graphicx}
\usepackage{xcolor}
\usepackage{tikz}
\usetikzlibrary{arrows.meta,positioning,shadows,calc}
\usepackage{caption}
\allowdisplaybreaks
\theoremstyle{plain}
\newtheorem{theorem}{Theorem}[section]
\newtheorem{proposition}[theorem]{Proposition}

\theoremstyle{definition}

\newtheorem{hypothesis}[theorem]{Hypothesis}
\newtheorem{example}[theorem]{Example}
\theoremstyle{remark}

\newcommand{\R}{\mathbb{R}}
\newcommand{\muP}{\mu^{+}} \newcommand{\muM}{\mu^{-}}
\newcommand{\nuP}{\nu^{+}} \newcommand{\nuM}{\nu^{-}}
\newcommand{\piPP}{\pi^{++}} \newcommand{\piMM}{\pi^{--}}
\newcommand{\piPM}{\pi^{+-}} \newcommand{\piMP}{\pi^{-+}}
\newcommand{\supp}{\mathrm{supp}}
\newcommand{\Lamt}{\widetilde\Lambda}

\title[Inter-sign optimal transport: synthesis]
{Regularity of the positional penalization function in inter-sign optimal transport on real measures \\[3mm]
\normalfont\itshape\small --- Synthesis ---}

\author{Bwo'nyahre Ba\"idi Barth\'el\'emy}
\address{Department of Mathematics and Computer Science, Faculty of Science, University of Ngaound\'er\'e, Cameroon}
\email{bwonyahre@proton.me}

\author{Kouakep Tchaptchie Yannick}
\address{School of Chemical Engineering and Mineral Industries (EGCIM), University of Ngaound\'er\'e, Cameroon}
\email{kouakep@aims-senegal.org}

\author{Houpa Danga Duplex Elvis}
\address{Department of Mathematics and Computer Science, Faculty of Science, University of Ngaound\'er\'e, Cameroon}
\email{e\_houpa@yahoo.com}


\begin{document}

\begin{abstract}
We announce a theory of optimal transport between \emph{signed} measures in which transport between opposite-sign masses is penalized by an additional, position-dependent cost $\lambda(x,y)$ on top of the ground cost $c(x,y)$. We state: existence and Kantorovich duality for this penalized problem; a canonical decomposition of optimal plans into intra-sign and inter-sign components; a variational characterization of the penalization function itself, $\lambda=\Lambda^*$; local Lipschitz and Alexandrov regularity of $\lambda$; and a modified Monge--Amp\`ere equation, in the form of a symmetric Sylvester equation, governing inter-sign transport maps, which reduces to the classical Brenier--Monge--Amp\`ere equation as $\lambda\to0$. All statements below are validated independently and quantitatively by a homemade entropic-transport (Sinkhorn) implementation. This is an announcement: proofs, together with a full numerical validation, appear in the companion article~\cite{full}.
\end{abstract}

\maketitle

\section{Motivation}

Classical optimal transport (Monge--Kantorovich) moves a nonnegative measure $\mu$ onto a nonnegative measure $\nu$ of equal mass at minimal cost $\int c\,d\Pi$. Extending this to \emph{signed} measures $\mu=\muP-\muM$, $\nu=\nuP-\nuM$ (Jordan decompositions) raises a structural question that has no classical analogue: a coupling may move mass from $\muP$ onto $\nuM$ (or from $\muM$ onto $\nuP$) -- an operation with no sign-preserving interpretation, amounting to a partial cancellation, or annihilation, of mass. We regulate this inter-sign transport by an explicit \emph{positional penalization} $\lambda(x,y)\ge0$, added to the ground cost precisely on inter-sign pairings. This yields a well-posed variational problem whose Kantorovich duality, canonical structure, and regularity theory are the object of this note.

\section{Setting}

Let $X,Y\subset\R^d$ and let $\mu=\muP-\muM$, $\nu=\nuP-\nuM$ be signed Radon measures on $X,Y$ with $\muP,\muM,\nuP,\nuM\ge0$ mutually singular and of matching total mass. Fix a ground cost $c$ and a penalization $\lambda\ge0$, and set $\Lamt:=c+\lambda$. An \emph{admissible plan} is a quadruple $\Pi=(\piPP,\piPM,\piMP,\piMM)$ of nonnegative measures on $X\times Y$ with
\[
  \mathrm{Proj}_X(\piPP+\piPM)=\muP,\quad
  \mathrm{Proj}_X(\piMM+\piMP)=\muM,\quad
  \mathrm{Proj}_Y(\piPP+\piMP)=\nuP,\quad
  \mathrm{Proj}_Y(\piMM+\piPM)=\nuM,
\]
and we write $\mathcal A(\mu,\nu)$ for the (nonempty) set of such plans. The cost is
\[
  \mathcal C(\Pi) = \int c\,d\piPP + \int c\,d\piMM + \int \Lamt\,d\piPM + \int \Lamt\,d\piMP,
\]
and we study $\inf_{\Pi\in\mathcal A(\mu,\nu)} \mathcal C(\Pi)$. Throughout, $X,Y$ are compact and $\muP,\muM,\nuP,\nuM$ are absolutely continuous with densities $f^\pm,g^\pm$; further structural hypotheses (convexity of $X,Y$, strict convexity and smoothness of $c$, $C^2$ densities) are added where stated, exactly as in~\cite{full}.

\section{Main results}

\subsection{Existence and duality}

\begin{theorem}[Existence]
$\mathcal A(\mu,\nu)\ne\emptyset$ and $\inf_{\mathcal A(\mu,\nu)}\mathcal C$ is attained.
\end{theorem}

\begin{theorem}[Kantorovich duality]\label{thm:duality}
\[
  \min_{\Pi\in\mathcal A(\mu,\nu)}\mathcal C(\Pi) =
  \sup\Bigl\{ \textstyle\int\varphi^+d\muP+\int\varphi^-d\muM+\int\psi^+d\nuP+\int\psi^-d\nuM \Bigr\}
\]
over Lipschitz potentials $(\varphi^+,\varphi^-,\psi^+,\psi^-)$ satisfying
\[
  \mathrm{(D1)}\ \varphi^++\psi^+\le c,\quad
  \mathrm{(D2)}\ \varphi^-+\psi^-\le c,\quad
  \mathrm{(D3)}\ \varphi^++\psi^-\le \Lamt,\quad
  \mathrm{(D4)}\ \varphi^-+\psi^+\le \Lamt.
\]
At optimality, complementary slackness holds on each of the four supports (e.g.\ $\piPM>0\Rightarrow\varphi^++\psi^-=\Lamt$).
\end{theorem}

\subsection{Canonical decomposition and the variational characterization of $\lambda$}

\begin{theorem}[Canonical decomposition]
Every optimal $\Pi$ admits a unique decomposition into intra-sign parts ($\piPP,\piMM$) and inter-sign parts ($\piPM,\piMP$), and this decomposition is stable under the marginal constraints.
\end{theorem}

Write $h^1:=\varphi^++\psi^--c$, $h^2:=\varphi^-+\psi^+-c$, and
$\Lambda^*:=\max(h^1,h^2)$.

\begin{theorem}[Variational characterization]\label{thm:regularity}
$\lambda\ge\Lambda^*$ on $X\times Y$, with equality on $\supp(\piPM)\cup\supp(\piMP)$: explicitly, $\lambda=h^1$ on $\supp(\piPM)$ and $\lambda=h^2$ on $\supp(\piMP)$.
\end{theorem}

\begin{proposition}[Compatibility]\label{prop:compat}
If $\supp(\piPM)\cap\supp(\piMP)\ne\emptyset$, then on this set,
$\varphi^+(x)-\varphi^-(x)=\psi^+(y)-\psi^-(y)$.
\end{proposition}

\subsection{Regularity of $\lambda$}

\begin{theorem}[Lipschitz and Alexandrov regularity]\label{thm:alex}
Under the structural hypotheses of~\cite{full} (in particular $X,Y$ convex compact, $c\in C^2$ strictly convex in $y$, $C^2$ strictly positive densities), $\lambda=\Lambda^*$ is locally Lipschitz on the interior of the inter-sign supports, and locally semiconvex there; consequently $\lambda$ admits classical second-order partial derivatives $D^2_{yy}\lambda$, $D^2_{xy}\lambda$ at Lebesgue-almost every point (Alexandrov).
\end{theorem}

An unconditional existence result for the inter-sign map (no ellipticity needed) underlies the next statement: under absolute continuity of $\muP$ and strict convexity of $\Lamt(x,\cdot)$, the optimal inter-sign plan $\piPM$ is induced by a Borel map $T:X^{+-}\to Y$ (Gangbo--McCann twist-condition argument). Ellipticity is then used only to upgrade $T$ to a $C^1$-diffeomorphism.

\subsection{A modified Monge--Amp\`ere equation}

\begin{hypothesis}[Ellipticity]\label{hyp:ell}
The matrix $A(x):=D^2_{yy}\Lamt(x,T(x))-D^2\psi^-(T(x))$ is uniformly positive definite on $X^{+-}$ (and analogously for $\piMP$). A local convexity condition on $\lambda$ in $y$ is sufficient.
\end{hypothesis}

\begin{theorem}[Modified Monge--Amp\`ere equation]\label{thm:MA}
Under Hypothesis~\ref{hyp:ell} and the Ma--Trudinger--Wang condition on $\Lamt$, $T$ is a $C^1$-diffeomorphism and satisfies, on $X^{+-}$,
\[
  \det\bigl(D^2_{yy}\Lamt(x,T(x))-D^2\psi^-(T(x))\bigr)
  = \frac{g^-(T(x))}{f^+(x)}\,\bigl|\det D^2_{xy}\Lamt(x,T(x))\bigr|,
\]
and analogously on $X^{-+}$. When $\lambda\equiv0$ this is the classical Brenier--Monge--Amp\`ere equation.
\end{theorem}

\begin{example}
For $c(x,y)=\tfrac12|x-y|^2$ and $\lambda(x,y)=\tfrac{\gamma}{2}(x-y)^2$ ($\gamma>0$, $d=1$), the equation reads $(1+\gamma)-\psi''(T(x)) = \dfrac{g^-(T(x))}{f^+(x)}(1+\gamma)$, manifestly elliptic since $f^+,g^->0$.
\end{example}

\section{Numerical validation}

\begin{figure}[h]
\centering
\begin{tikzpicture}[
    >=Stealth, line width=0.8pt, scale=0.82, every node/.style={scale=0.82},
    box/.style={draw, rounded corners=3pt, minimum width=2.5cm, minimum height=1.0cm,
                font=\small, align=center, drop shadow={opacity=0.15,shadow xshift=1pt,shadow yshift=-1pt}},
    lab/.style={font=\scriptsize, midway, above, sloped}
]
\definecolor{navy}{HTML}{0b3d66}
\definecolor{red1}{HTML}{a3123a}
\definecolor{green1}{HTML}{1a7a4c}
\definecolor{gold1}{HTML}{b8860b}
\node[box, fill=navy!8, draw=navy] (Xp) at (0,1.9)  {$X^+=X\times\{+\}$\\[1pt]{\scriptsize mass $\muP$}};
\node[box, fill=green1!8, draw=green1] (Xm) at (0,-1.9) {$X^-=X\times\{-\}$\\[1pt]{\scriptsize mass $\muM$}};
\node[box, fill=navy!8, draw=navy] (Yp) at (5.7,1.9)  {$Y^+=Y\times\{+\}$\\[1pt]{\scriptsize mass $\nuP$}};
\node[box, fill=green1!8, draw=green1] (Ym) at (5.7,-1.9) {$Y^-=Y\times\{-\}$\\[1pt]{\scriptsize mass $\nuM$}};
\draw[->, navy]  (Xp) -- (Yp) node[lab, navy]   {cost $c$\; ($\piPP$)};
\draw[->, green1] (Xm) -- (Ym) node[lab, green1] {cost $c$\; ($\piMM$)};
\draw[->, red1]  (Xp) to[bend left=10] (Ym);
\draw[->, gold1] (Xm) to[bend right=10] (Yp);
\node[red1, font=\scriptsize, anchor=north east] at ($(Ym.north west)+(-0.05,0.35)$) {$\piPM$};
\node[gold1, font=\scriptsize, anchor=south east] at ($(Yp.south west)+(-0.05,-0.35)$) {$\piMP$};
\begin{scope}[shift={(9.0,-1.05)}]
  \draw[navy, fill=navy!12] (0,1) rectangle (1,2);
  \draw[green1, fill=green1!12] (1,0) rectangle (2,1);
  \draw[red1, fill=red1!12] (1,1) rectangle (2,2);
  \draw[gold1, fill=gold1!12] (0,0) rectangle (1,1);
  \node[font=\scriptsize] at (0.5,1.5) {$c$};
  \node[font=\scriptsize] at (1.5,0.5) {$c$};
  \node[font=\scriptsize] at (1.5,1.5) {$c+\lambda$};
  \node[font=\scriptsize] at (0.5,0.5) {$c+\lambda$};
  \node[font=\scriptsize] at (-0.55,1.5) {$X^+$};
  \node[font=\scriptsize] at (-0.55,0.5) {$X^-$};
  \node[font=\scriptsize] at (0.5,-0.35) {$Y^+$};
  \node[font=\scriptsize] at (1.5,-0.35) {$Y^-$};
  \node[font=\scriptsize] at (1.0,2.35) {block cost $\widehat C$};
\end{scope}
\end{tikzpicture}
\caption{Doubled-space reformulation used by the homemade Sinkhorn solver: straight arrows are intra-sign (cost $c$), curved arrows inter-sign (cost $c+\lambda$), color-matched to the block cost $\widehat C$ on the right.}
\label{fig:doubled}
\end{figure}
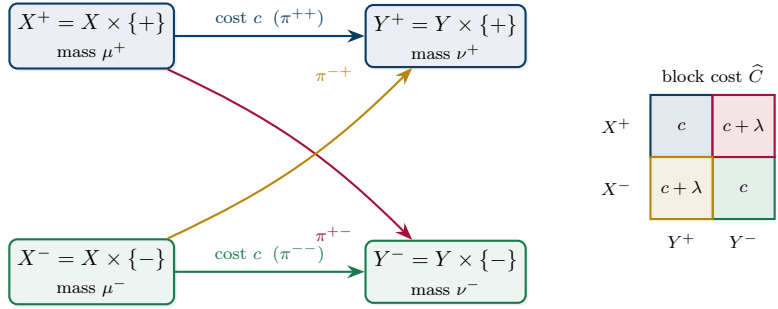

Every statement above is checked independently and quantitatively in~\cite{full} through a homemade log-domain Sinkhorn solver (no external optimal-transport library), run on a doubled-space reformulation $\widehat X=X\times\{+,-\}$, $\widehat Y=Y\times\{+,-\}$ that turns the four-block problem into a single standard entropic transport problem (Figure~\ref{fig:doubled}): a standard coupling $\widehat\pi$ between $\widehat\mu=\muP\oplus\muM$ and $\widehat\nu=\nuP\oplus\nuM$, for the block cost shown below, recovers exactly the four blocks $(\piPP,\piPM,\piMP,\piMM)$ of Section~2.

\begin{figure}[h]
\centering
\includegraphics[width=0.62\linewidth]{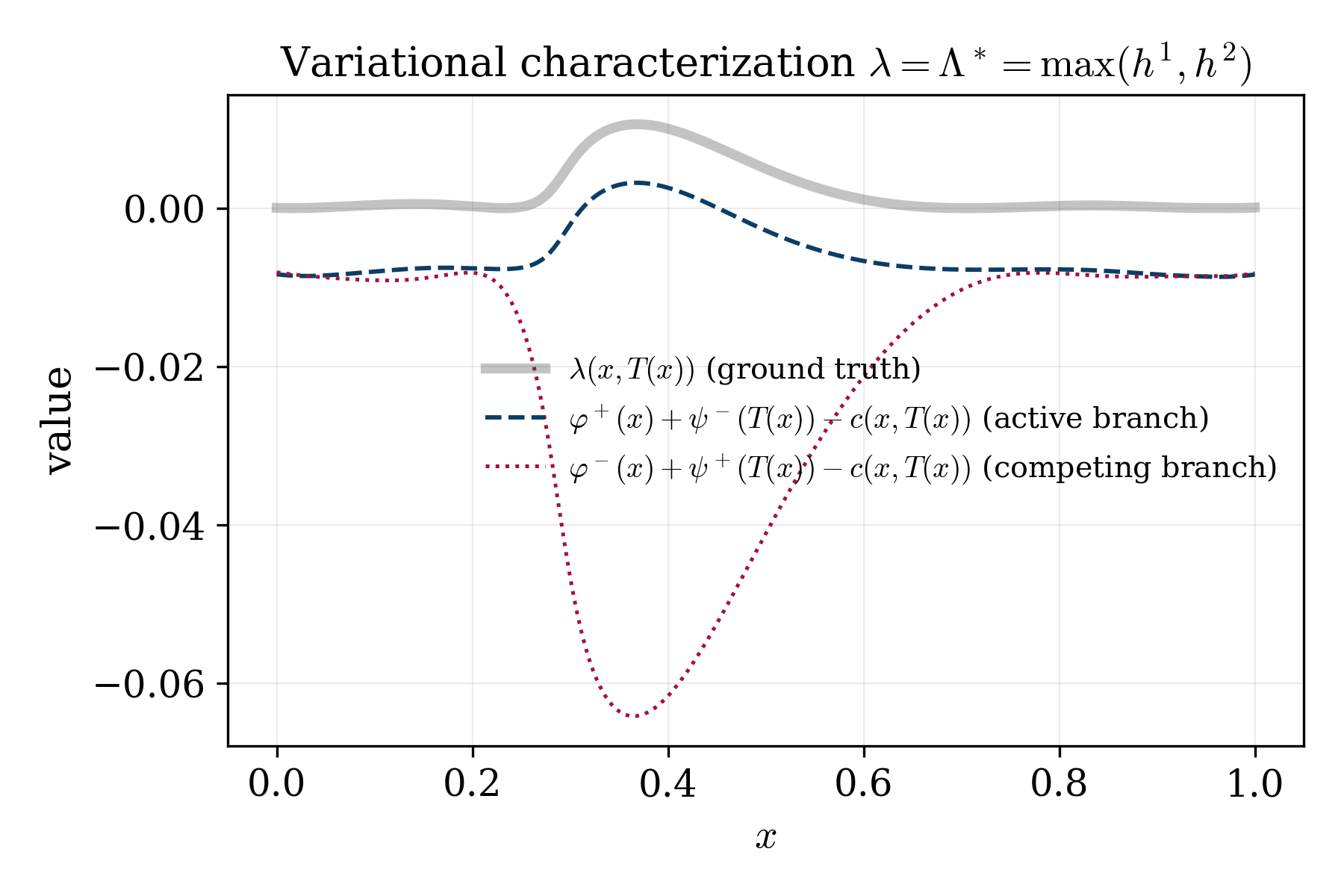}
\caption{$d=1$: $\lambda(x,T(x))$ (ground truth) versus the active branch $h^1=\varphi^++\psi^--c$ and the competing branch $h^2$, recovered from the Sinkhorn potentials, confirming Theorem~\ref{thm:regularity}.}
\label{fig:lambdastar}
\end{figure}

Ground truths are exact: the monotone rearrangement in $d=1$, and Gaussian-to-Gaussian transport (whose Brenier map is affine, and independent of positive rescalings of the cost) in $d=2$, which gives closed-form, matrix-valued Hessians for a genuine test of Theorem~\ref{thm:MA}. Figure~\ref{fig:lambdastar} illustrates the central check: the variational characterization $\lambda=\Lambda^*=\max(h^1,h^2)$ of Theorem~\ref{thm:regularity}, recovered from the Sinkhorn potentials alone. Table~\ref{tab:residuals} summarizes all residuals.

\begin{table}[h]
\centering
\caption{Numerical residuals (median absolute value, $\varepsilon\to0$).}
\label{tab:residuals}
\begin{tabular}{@{}lc@{}}
\toprule
Quantity checked & Residual \\
\midrule
$\lambda=\Lambda^*$ (Theorem~\ref{thm:regularity}), $d=1$ & $7.8\times10^{-3}$ \\
Compatibility (Proposition~\ref{prop:compat}), $d=1$ & $5.8\times10^{-4}$ \\
Modified Monge--Amp\`ere equation, $d=1$ & $1.1\times10^{-2}$ \\
Modified Monge--Amp\`ere equation (matrix, $\det$), $d=2$ Gaussian & $3.2\times10^{-2}$ (relative) \\
\bottomrule
\end{tabular}
\end{table}

\section{Perspectives}

The absolutely continuous case treated here and in~\cite{full} is a first step; the general case, allowing $\lambda$-penalized transport between signed measures with fractal singular parts, is the subject of a companion article~\cite{fractal}. On the numerical side, extending the quantitative $d\ge2$ benchmark beyond the Gaussian case, and a direct numerical check of the Ma--Trudinger--Wang condition itself, are left for future work.

\end{document}